# A HYBRID COA-DEA METHOD FOR SOLVING MULTI-OBJECTIVE PROBLEMS

Mahdi Gorjestani[1], Elham shadkam[2], Mehdi Parvizi[3] and Sajedeh Aminzadegan[4]

[1,3,4]Department of Industrial Engineering, Faculty of Eng., Khayyam University, Mashhad, Iran
[2] Ph.D. Candidate of Department of Industrial and Systems Engineering, Isfahan University of Technology, Isfahan, Iran, and Faculty Member of Industrial Engineering Department, Faculty of Eng, Khayyam University, Mashhad, Iran

## Abstract

*The Cuckoo optimization algorithm (COA) is developed for solving single-objective problems and it cannot be used for solving multi-objective problems. So the multi-objective cuckoo optimization algorithm based on data envelopment analysis (DEA) is developed in this paper and it can gain the efficient Pareto frontiers. This algorithm is presented by the CCR model of DEA and the output-oriented approach of it. The selection criterion is higher efficiency for next iteration of the proposed hybrid method. So the profit function of the COA is replaced by the efficiency value that is obtained from DEA. This algorithm is compared with other methods using some test problems. The results shows using COA and DEA approach for solving multi-objective problems increases the speed and the accuracy of the generated solutions.*

## 1.Introduction

Finding the best solution for an objective subject to some conditions calls optimization. In multi-objective problems, there is not an optimal solution that can optimize all objectives simultaneously. So, in order to solve problems the concept of Pareto frontiers is provided. Usually, there are some Pareto optimized solutions that the best solution will be selected from them by decision maker. Many practical problems in real world are multi-objective problem. Several researches developed for solving multi-objective problem.

Ehrgott and Gandibleux studied on the approximate and the accurate problems related to the combination method of multi-objective problems [1]. Arakawa *et al.* combined the General DEA and the Genetic Algorithm to generate the efficient frontier in multi-objective optimization problems [2]. Deb used the evolutionary algorithms for solving the multi-objective problem [3].

Nakayama et al. drew the Pareto frontier of the multi-objective optimization problems using DEA in 2001 [4]. Deb et al. obtained the Pareto frontier of the multi-objective optimization problem using Genetic Algorithm [5]. Kristina Vincova gained the Pareto frontier using DEA [6]. Reyes-Sierra and Coello Coello investigated the method to solve the multi-objective optimization using the particle swarm algorithm [7]. Cooper et al. and Tone improved the multi-objective optimization algorithm using DEA and developed related software [8]. Pham and Ghanbarzadeh solved the multi-objective optimization algorithm using the Bees Algorithm [9]. Nebro et al. investigated a new method for multi-objective optimization algorithm based on the particle swarm algorithm [10]. Yun *et al.* studied the solution of multi-objective optimization algorithm using the GA and DEA. Also, they applied their method to generating the Pareto efficient frontiers [11]. Yang *and* Deb used the Cuckoo optimization algorithm in order to solve the multi objective-problem [12].

In this article, it is tried to use the meta-heuristic Cuckoo algorithm with the DEA approach for solving multi-objective problems and draw the Pareto frontiers for efficient points of the considered objective functions. Because of using CCR model of DEA, proposed method only applicable to generating the convex efficient frontier. In the second section, the Cuckoo algorithm is introduced. In the third section the multi-objective problems are defined. In the fourth section, the concept of DEA is explained. The fifth section expresses the proposed hybrid method and in the sixth section the test problems are given. At last the desired conclusion is provided.

## 2.Introducing the Cuckoo optimization algorithm

The COA is one of the best and newest evolutionary algorithms. After early evolutionary methods like Genetic algorithm (GA), Simulated Annealing algorithm, so many evolutionary methods that inspired from the nature, have been developed. Some of the useful algorithms for solving complicated optimization problems are Particle Swarm Optimization (PSO), Ant Colony Optimization (ACO), Artificial Bee Colony algorithm (ABC) and the Artificial Fish Swarm algorithm. One of the other evolutionary algorithms that are developed in Iran is Imperialist Competitive Algorithm (ICA). This algorithm is inspired from the competitive system of the empires in order to get more colonies. After the ICA, the Cuckoo optimization algorithm is presented that has the ability to find the general optimized solutions. This algorithm is inspired from the life of a bird calls Cuckoo. The cuckoo living and egg laying method is a suitable inspiration for inventing an evolutionary algorithm. The survival with the least effort is the base of this method. This lazy bird forces other birds to play an important role in her survival so nicely. The Cuckoo optimization algorithm expanded by Yang and Deb in 2009. This algorithm is inspired by the egg laying method of cuckoos combined with the Levy Flight instead of simple random isotropic walk. The COA investigated with more details by Rajabioun in 2011 [13].
The flowchart of the COA is given in the Figure 1.

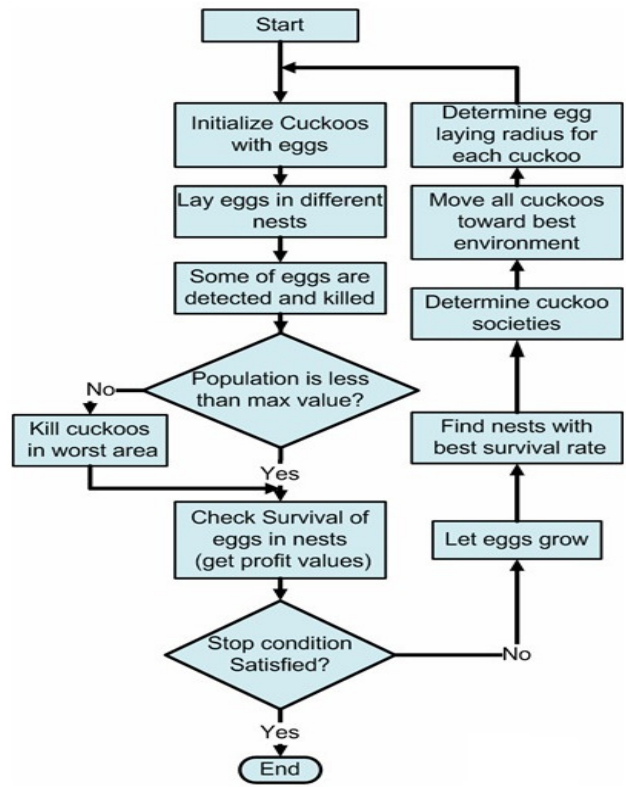

Figure1. The flowchart of COA algorithm

For more information refers to [13].

## 3.The multi-objective optimization problem

The general form of a multi-objective optimization problem is as (1):

$$Max \ or \ Min = \{f_1, f_2, \dots, f_k\}$$
$$s.t. \ \ g(x) = \sum_{j=1}^{n} c_j x_j \leq b_i, \qquad i = 1, \dots, m$$
$$x_j \geq 0, \qquad j = 1, \dots, n$$

(1)

As it shown, we encounter to several numbers of the objectives in multi-objective problem. $K$ is the number of objective functions that can be min or max type, $m$ is the number of constraints and $n$ is the number of decision variables. In multi-objective algorithms, there is not an optimal solution that can optimize all of the objective functions simultaneously. So the concept of Pareto optimized solution is provided. The Pareto optimal concept is explicable this way. $\overline{x_*} = (x_1 \dots$ ,$x_2, x_n)$ is an optimal Pareto, if for each allowable x̄ and i={1,2,..k}, we have (for minimizing problem is as (2)):

$$\forall_{i \in I}(f_i(\overline{x}_*) \leq f_i(\overline{x}_i)) \tag{2}$$

In other words, $\overline{x}_*$ is an optimal Pareto solution when no other $\overline{x}$ exist that make at least one objective function in order to optimizing some of the objective functions.

## 4.Data Envelopment Analysis

DEA is a linear programming technique. Its main purpose is comparing and evaluating the efficiency of a number of similar decision making units like banks, hospitals, schools, refineries, power plants etc. that the amount of their consumed input and production output are different.

DEA shows the concept of evaluating the efficiency within a group of DMUs. In this method the efficiency of each DMU is calculated according to other DMUs that have the most operations.

The first model of DEA is CCR, defining the efficiency according to the ratio of an output on an input, is the base of this method [6]. In other words, calculating the ratio of total weighed of outputs on total weighed of inputs instead of the ratio of an output on an input is used for evaluating the efficiency in CCR model.

### The CCR model

The CCR model is the first model of DEA and its named is the first letters of the names of its providers (Charnes, Cooper, Rhodez) [4]. For determining the best efficient unit, the amounts of inputs and outputs of other decision making units in finding the optimal weights for each unit is evaluated. This basic model is suggested as (3):

$$
\begin{aligned}
&Min \sum_{i=1} V_i x_{io} \\
&\sum_{r=1}^{k} u_r y_{rj} - \sum_{i=1}^{m} v_i x_{ij} \leq 0 \qquad j=1,...,n \\
&s.t. \sum_{r=1}^{k} u_r y_{ro} = 1 \\
&u_r \geq 0 \qquad v_i \geq 0
\end{aligned}
\tag{3}
$$

Where $u_r$ is the weight of output $r$, $v_i$ is the weight of input $i$ and $o$ is the index of under reviewed DMU, $(o \in \{1,2,...,n\})$. $y_{ro}$ is the amount of output $r$ and $x_{io}$ is the amount of input $i$ for DMU$_o$. Also $y_{ij}$ is the amount of output $r$ and $x_{ij}$ is the amount of input $i$ for the unit $j$. $k$ is the number of outputs; $m$ is the number of inputs and $n$ is the number of units.

# 5. The proposed hybrid algorithm

In this paper, it is tried to present a hybrid method in order to solve the multi-objective problems using COA and DEA methods. This hybrid method finds the efficient points using DEA method and gains the Pareto frontiers for multi-objective problems.

## The steps of hybrid COA_DEA algorithm

1. In the first step of implementing the Cuckoo algorithm, the desired matrix will be formed from habitats according to the initial population of cuckoos and the initial egg laying radiuses.
2. The "profit function" of the Cuckoo algorithm will be replaced by the "efficiency value". This function take the habitat matrix as its input according to this matrix, the CCR model will be produced for each habitats of the matrix and determines the efficiency for each habitat.
3. The habitats will be sorted according to their efficiency values and other steps will be as the explanations that are given in the references [13].
4. In each iteration, the habitats with the efficiency of one will be selected as good solutions for transferring to next iteration.
5. At last iteration of the proposed algorithm, The Pareto frontiers for the main multi-objective optimization problem will be drawn out based on the obtained values of $f_1$ and $f_2$.

# 6. SOLVING TEST PROBLEMS

A number of test functions have been provided that can help to validate the proposed method in Table 1.

Table 1. Test problems

| Number | Objective function | Constraints | Figure number |
|---|---|---|---|
| 1 | $minf_1 = 4x_1 + 4x_2$ <br> $minf_2 = (x_1 - 5)^2 + (x_2 - 5)^2$ | $(x_1 - 5)^2 + x_2{}^2 - 25 \leq 0$ | (1) |
| 2 | $minf_1 = 2x_1 - x_2$ <br> $minf_2 = -x_1$ | $(x_1 - 1)^3 + x_2 \leq 0 \, ;$ | (2) |
| 3 | $min f1 = (x_1 - 2)^2 + 2 + (x_2 - 1)^2 + 2$ <br> $min f_2 = 9x_1 + (x_2 - 1)^2$ | $x_1{}^2 + x_2{}^2 \leq 225$ <br> $x_1 + 3x_2 \leq -10$ | (3) |
| 4 | $min f_1 = x_1$ <br> $min f_2 = \dfrac{1 + x_2}{x_1}$ | $x_2 + 9x_1 \geq 6$ <br> $-x_2 + 9x_1 \geq 1$ | (4) |

Parameters setting for cuckoo algorithm are as follow:

Number of initial population=5, minimum number of eggs for each cuckoo=2, maximum number of eggs for each cuckoo=6, maximum iterations of the Cuckoo Algorithm=8, number of clusters that we want to make=2, maximum number of cuckoos that can live at the same time=50.

***Test problem 1: [16]***

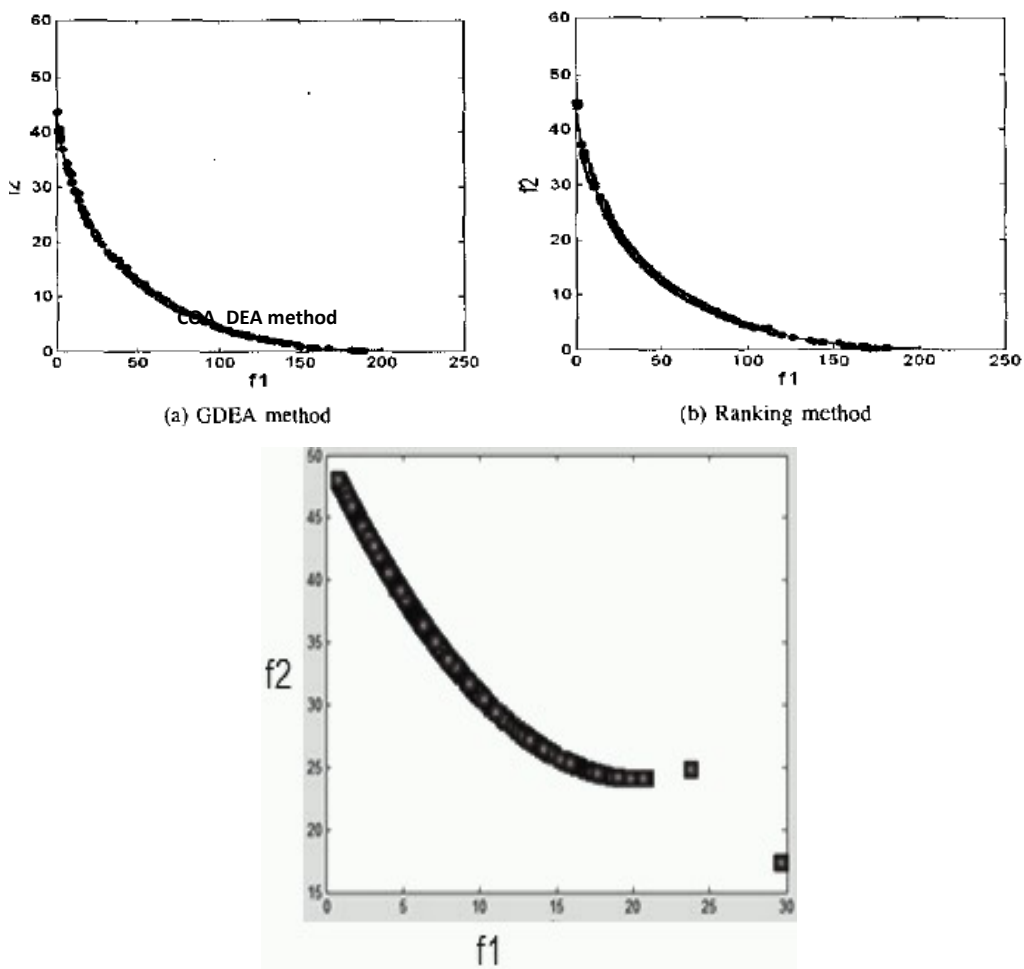

Figure 1. Comparing the proposed method with other methods

***Test problem 2: [14]***

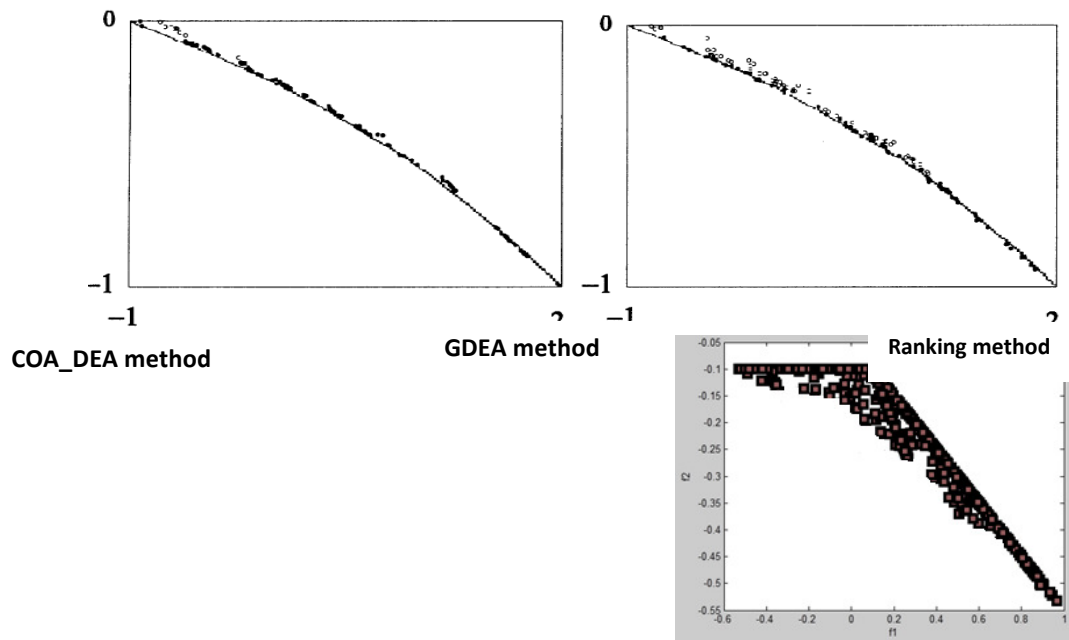

**COA_DEA method**        **GDEA method**        **Ranking method**

Figure2. Comparing the proposed method with other methods

***Test problem 3: [15]***

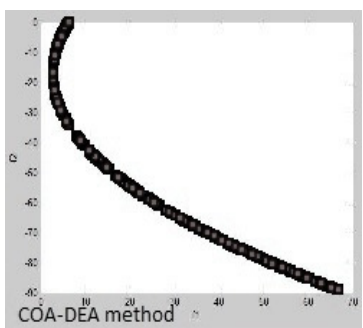

COA-DEA method

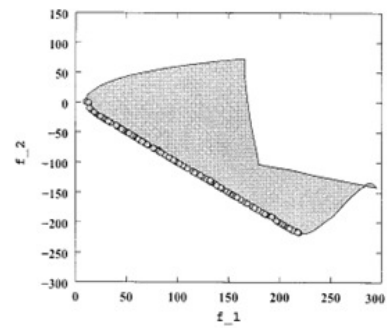

Fig. 16. Obtained nondominated solutions with NSGA-II on the constrained problem SRN.

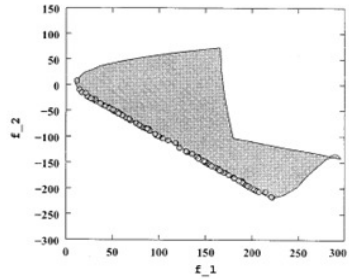

Fig. 17. Obtained nondominated solutions with Ray-Tai-Seow's algorithm on the constrained problem SRN.

Figure 3.Comparing the proposed method with other methods

*Test problem 4: [15]*

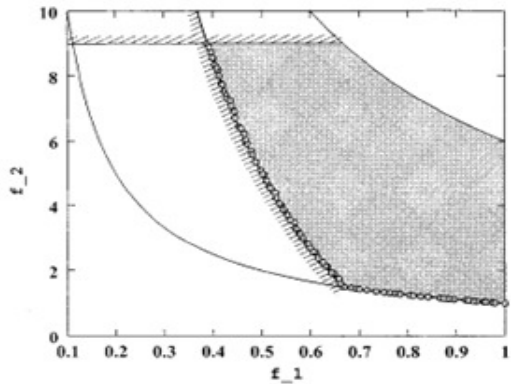

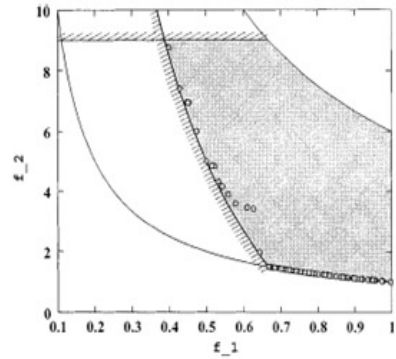

Fig. 14. Obtained nondominated solutions with NSGA-II on the constrained problem CONSTR.

Fig. 15. Obtained nondominated solutions with Ray-Tai-Seow's algorithm on the constrained problem CONSTR.

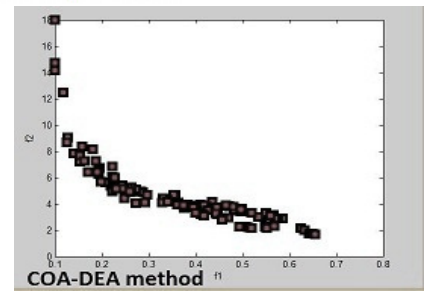

COA-DEA method

Figure 4. Comparing the proposed method with other methods

# 7.Conclusion

In this paper, it is tried to solve multi-objective problems with a new creative approach. This approach is a combination of the Cuckoo optimization algorithm and DEA method. As it shown this method is one of the fastest, most accurate and most logical method for solving multi-objective problems because it is a logical combination of both efficiency and finding the optimal solutions. We conclude that the proposed method not only finds optimal answers and more efficient points, but also it is faster in processing time than other algorithms. The obtained Pareto frontiers of this method were compared with the answers of similar algorithms like GA-DEA, Ranking method, GA-GDEA, etc. The algorithm's convergence rate in order to find the answer is evident. So the suggested method is suitable and reliable method for solving multi-objective optimization problems.

For further work, we can use another clustering method instead of current method for grouping the cuckoos.